\input gtmacros
\input amsnames
\input amstex
%
\catcode`\@=12        
%
%
\input gtoutput
\volumenumber{3}\papernumber{3}\volumeyear{1999}
\pagenumbers{67}{101}\published{28 May 1999}
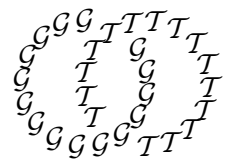
\proposed{Ralph Cohen}\seconded{Haynes Miller, Gunnar Carlsson}
\received{10 May 1998}\revised{5 May 1999}
\accepted{13 May 1999}
%
\let\\\par
\def\topmatter{\relax}
\def\endtopmatter{\maketitlepage}
\let\gttitle\title
\def\title#1\endtitle{\gttitle{#1}}
\let\gtauthor\author
\def\author#1\endauthor{\gtauthor{#1}}
\let\gtaddress\address
\def\address#1\endaddress{\gtaddress{#1}}
\let\gtemail\email
\def\email#1\endemail{\gtemail{#1}}
\def\subjclass#1\endsubjclass{\primaryclass{#1}}
\let\gtkeywords\keywords
\def\keywords#1\endkeywords{\gtkeywords{#1}}
\def\heading#1\endheading{{\def\S##1{\relax}\def\\{\relax\ignorespaces}
    \section{#1}}}
\def\head#1\endhead{\heading#1\endheading}
\def\subheading#1{\sh{#1}}
\def\subhead#1\endsubhead{\sh{#1}}
\def\subsubhead#1\endsubsubhead{\sh{#1}}
\def\specialhead#1\endspecialhead{\sh{#1}}
\def\demo#1{\rk{#1}\ignorespaces}
\def\enddemo{\ppar}
\let\remark\demo
\def\endremark{}
\let\definition\demo
\def\enddefinition{\ppar}

\def\qed{\ifmmode\quad\sq\else\hbox{}\hfill$\sq$\par\goodbreak\rm\fi}  
\def\proclaim#1{\rk{#1}\sl\ignorespaces}
\def\endproclaim{\rm\ppar}
\def\cite#1{[#1]}
\newcount\itemnumber
\def\roster{\items\itemnumber=1}
\def\endroster{\enditems}
\def\therosteritem#1{{\rm(#1)}}
\let\itemold\item
\def\item{\itemold{{\rm(\number\itemnumber)}}%
\global\advance\itemnumber by 1\ignorespaces}
\def\S{section~\ignorespaces}  
\def\date#1\enddate{\relax}
\def\thanks#1\endthanks{\relax}   
\def\dedicatory#1\enddedicatory{\relax}  
\let\footnote\plainfootnote

\def\Refs{\ppar{\large\bf References}\ppar\bgroup\leftskip=25pt
\frenchspacing\parskip=3pt plus2pt\small}       
\def\endRefs{\egroup}
\def\widestnumber#1#2{\relax}
\def\endrefitem{}
\def\refdef#1#2#3{\def#1{\leavevmode\unskip\endrefitem#2\def\endrefitem{#3}}} 
\def\ref{\par}
\def\endref{\endrefitem\par\def\endrefitem{}}
\refdef\key{\noindent\llap\bgroup[}{]\ \ \egroup}
\refdef\no{\noindent\llap\bgroup[}{]\ \ \egroup}
\refdef\by{\bf}{\rm, }
\refdef\manyby{\bf}{\rm, }
\refdef\paper{\it}{\rm, }
\refdef\book{\it}{\rm, }
\refdef\jour{}{ }
\refdef\vol{}{ }
\refdef\yr{$(}{)$ }
\refdef\ed{(}{ Editor) }
\refdef\publ{}{ }
\refdef\inbook{from: ``}{'', }
\refdef\pages{}{ }
\refdef\page{}{ }
\refdef\paperinfo{}{ }
\refdef\bookinfo{}{ }
\refdef\publaddr{}{ }
\refdef\eds{(}{ Editors)}
\refdef\bysame{\hbox to 3 em{\hrulefill}\thinspace,}{ }
\refdef\toappear{(to appear)}{ }
\refdef\issue{no.\ }{ }
\refdef\moreref{ }{ }

\define\la{\longrightarrow}
\define\bla{\text{---}}
\define\shtimes{\!\times\!}
\define\minus{\smallsetminus}
\define\down{\!\downarrow\!}

\define\cO{\Cal O}
\define\cF{\Cal F}
\define\cP{\Cal P}
\define\cX{\Cal X}

\define\cT{\Cal T}
\define\cC{\Cal C}
\define\cA{\Cal A}
\define\cB{\Cal B}
\define\cI{\Cal I}
\define\cQ{\Cal Q}
\define\cJ{\Cal J}
\define\cR{\Cal R}
\define\cY{\Cal Y}

\define\RR{\Bbb R}
\define\DD{\Bbb D}

\define\spr{\operatorname{sp}}
\define\id{\operatorname{id}}
\define\im{\operatorname{im}}
\define\mor{\operatorname{mor}}
\define\nat{\operatorname{nat}}
\define\holim{\operatornamewithlimits{holim}}
\define\colim{\operatornamewithlimits{colim}}
\define\hocolim{\operatornamewithlimits{hocolim}}
\define\hofiber{\operatorname{hofiber}}
\define\imm{\operatorname{imm}}
\define\emb{\operatorname{emb}}
\define\map{\operatorname{map}}
\define\mono{\operatorname{mono}}

\define\intr{\operatorname{int}}

\topmatter
\title Embeddings from the point of view of\\immersion theory : Part I 
\endtitle
\shorttitle{Embeddings from immersion theory : I}
\author 
Michael Weiss
\endauthor
\address 
Department of Mathematics, University of Aberdeen\\Aberdeen, AB24 3UE, UK
\endaddress
\email m.weiss\@maths.abdn.ac.uk \endemail
\primaryclass{57R40}\secondaryclass{57R42}
\keywords Embedding, immersion, calculus of functors
\endkeywords
\abstract
Let $M$ and $N$ be smooth manifolds without 
boundary. Immersion theory suggests that an understanding 
of the space of smooth embeddings $\emb(M,N)$ 
should come from an analysis of the cofunctor
$V\mapsto\emb(V,N)$ from the poset $\cO$ of open subsets of $M$ to spaces.
We therefore abstract some of the properties of this cofunctor, and 
develop a suitable calculus of such cofunctors, Goodwillie style, with Taylor 
series and so on. The terms of the Taylor series for the cofunctor 
$V\mapsto\emb(V,N)$ are explicitly determined. In a sequel to this 
paper, we introduce the concept of an analytic cofunctor from $\cO$ 
to spaces, and show that the Taylor series of an analytic cofunctor $F$
converges to $F$. Deep excision theorems due to Goodwillie and 
Goodwillie--Klein imply that the cofunctor $V\mapsto\emb(V,N)$
is analytic when $\dim(N)-\dim(M)\ge3$. 
\endabstract
\asciiabstract{Let M and N be smooth manifolds without 
boundary. Immersion theory suggests that an understanding 
of the space of smooth embeddings emb(M,N) 
should come from an analysis of the cofunctor
V |--> emb(V,N) from the poset O of open subsets of M to spaces.
We therefore abstract some of the properties of this cofunctor, and 
develop a suitable calculus of such cofunctors, Goodwillie style, with Taylor 
series and so on. The terms of the Taylor series for the cofunctor 
V |--> emb(V,N) are explicitly determined. In a sequel to this 
paper, we introduce the concept of an analytic cofunctor from O 
to spaces, and show that the Taylor series of an analytic cofunctor F
converges to F. Deep excision theorems due to Goodwillie and 
Goodwillie-Klein imply that the cofunctor V |--> emb(V,N) 
is analytic when dim(N)-dim(M) > 2.}
\endtopmatter

\newpage 

\document
\catcode`\@=\active

\sectionnumber=-1
\section{Introduction} 
Recently Goodwillie \cite{9}, \cite{10}, \cite{11} and 
Goodwillie--Klein \cite{12}
proved higher excision theorems of Blakers--Massey type 
for spaces of smooth embeddings.
In conjunction with a calculus framework,
these lead to a calculation of such spaces when
the codimension is at least 3. Here the goal is to set up the 
calculus framework. This is very similar to Goodwillie's
calculus of homotopy functors \cite{6}, \cite{7}, \cite{8},
but it is not a special case. Much of it has been known
to Goodwillie for a long time. For some history and a slow 
introduction, see \cite{23}. If a reckless introduction is required,
read on---but be prepared for Grothendieck topologies 
\cite{18} and homotopy limits \cite{1}, \cite{23, \S1}. 

Let $M$ and $N$ be smooth manifolds without boundary. Write $\imm(M,N)$
for the space of smooth immersions from $M$ to $N$. Let 
$\cO$ be the poset of open subsets of $M$, ordered by inclusion. 
One of the basic ideas
of immersion theory since Gromov \cite{14}, \cite{16}, \cite{19} 
is that $\imm(M,N)$ should be regarded 
as just one value of the cofunctor $V\mapsto\imm(V,N)$ from $\cO$ to
spaces. Here $\cO$ is treated as a category, with exactly one morphism
$V\to W$ if $V\subset W$, and no morphism if $V\not\subset W$; a 
{\sl cofunctor} is a contravariant functor. 

The poset or category $\cO$ has a Grothendieck topology
\cite{18, III.2.2} which we denote by $\cJ_1$. Namely, a family of 
morphisms $\{V_i\to W\mid i\in S\}$ qualifies as a {\sl covering} in $\cJ_1$ 
if every point of $W$ is contained in some $V_i$. More generally,
for each $k>0$ there is a Grothendieck topology $\cJ_k$ on $\cO$ in which
a family of morphisms $\{V_i\to W\mid i\in S\}$ qualifies as a {\sl covering}
if every finite subset of $W$ with at most $k$ elements is contained in 
some $V_i$. We will say that a 
cofunctor $F$ from $\cO$ to spaces is a {\sl homotopy sheaf} with respect 
to the Grothendieck topology $\cJ_k$ if for any covering
$\{V_i\to W\mid i\in S\}$
in $\cJ_k$ the canonical map 
$$F(W)\la \holim_{\emptyset\ne R\subset S}F(\cap_{i\in R}V_i)$$ 
is a homotopy equivalence. Here $R$ runs through the 
finite nonempty subsets of $S$. In view of the homotopy invariance
properties of homotopy inverse limits, the condition means that the
values of $F$ on large open sets are 
sufficiently determined for the homotopy theorist by the behavior of
$F$ on certain small open sets; however, it depends on $k$ how
much smallness we can afford.--- The main theorem of immersion
theory is that the cofunctor $V\mapsto\imm(V,N)$ from $\cO$ to spaces
is a homotopy sheaf with respect to $\cJ_1$, provided 
$\dim(N)$ is greater than $\dim(M)$ or
$\dim(M)=\dim(N)$ and $M$ has no compact component.

In this form, the theorem may not be very recognizable. It can be 
decoded as follows. Let $Z$ be the space of all 
triples $(x,y,f)$ where $x\in M$, $y\in N$ and $f\co T_xM\to T_yN$ is a 
linear monomorphism. Let $p\co Z\to M$ be the projection to the
first coordinate. For $V\in\cO$ we denote by $\Gamma(p~;V)$ the space 
of partial sections of $p$ defined over $V$. It is not hard to 
see that $V\mapsto\Gamma(p~;V)$ is a homotopy sheaf with 
respect to $\cJ_1$. (Briefly:  if
$\{V_i\to W\}$ is a covering in $\cJ_1$, then the canonical map
$q\co \hocolim_R\cap_{i\in R}V_i\to W$
is a homotopy equivalence according to \cite{24}, so that 
$\Gamma(p~;W)\,\simeq\,\Gamma(q^*p)\,\cong\,
\holim_R\Gamma(p~;\cap_{i\in R}V_i)\,.\,)$ There is an obvious inclusion
$$\imm(V,N)\hookrightarrow \Gamma(p~;V)\tag{$*$}$$
which we regard as a natural transformation between cofunctors in the
variable $V$. We want to show that \thetag{$*$} is a homotopy equivalence for every
$V$, in particular for $V=M$; this is the decoded version
of the main theorem of immersion 
theory, as stated in Haefliger--Poenaru \cite{15} for example 
(in the PL setting).
By inspection, \thetag{$*$} is indeed a 
homotopy equivalence when $V$ is diffeomorphic to $\RR^m$. An 
arbitrary $V$ has a smooth triangulation and can then be covered by the 
open stars $V_i$ of the triangulation. Since \thetag{$*$} is a homotopy 
equivalence for the $V_i$ and their finite intersections, it is a 
homotopy equivalence for $V$ by the homotopy sheaf property.   

Let us now take a look at the space of smooth embeddings
$\emb(M,N)$ from the same point of view. As before, we think
of $\emb(M,N)$ as just one value of the cofunctor $V\mapsto\emb(V,N)$
from $\cO$ to spaces. The cofunctor is clearly not a homotopy
sheaf with respect to the Grothendieck topology $\cJ_1$, except in 
some very trivial cases. For if it were, the inclusion
$$\emb(V,N)@>\subset>>\imm(V,N)\tag{$**$}$$
would have to be a homotopy equivalence for every $V\in\cO$, since it is
clearly a homotopy equivalence when $V$ is diffeomorphic to $\RR^m$. 
In fact it is quite appropriate to think of the cofunctor $V\mapsto\imm(V,N)$
as the {\sl homotopy sheafification} of $V\mapsto\emb(V,N)$, again
with respect to $\cJ_1$. The natural transformation \thetag{$**$} has a suitable
universal property which justifies the terminology.

Clearly now is the time to try out the smaller Grothendieck topologies $\cJ_k$
on $\cO$. For each $k>0$ the cofunctor $V\mapsto\emb(V,N)$ has
a homotopy sheafification with respect to $\cJ_k$.
Denote this by $V\mapsto T_k\emb(V,N)$. Thus $V\mapsto T_k\emb(V,N)$ is 
a homotopy sheaf on $\cO$ with respect to $\cJ_k$ 
and there is a natural transformation
$$\emb(V,N)\la T_k\emb(V,N)\tag{$*\!*\!*$}$$
which should be regarded as the {\sl best approximation} of 
$V\mapsto\emb(V,N)$ by a cofunctor which is a homotopy sheaf with respect 
to $\cJ_k$. I do not know of any convincing 
geometric interpretations of $T_k\emb(V,N)$ except of course in the case 
$k=1$, which we have already discussed. As Goodwillie 
explained to me, his excision theorem for diffeomorphisms
\cite{9}, \cite{10}, \cite{11} 
and improvements due to Goodwillie--Klein \cite{12} imply that \thetag{$*\!*\!*$} is 
$(k(n-m-2)+1-m)$--connected where $m=\dim(M)$ and $n=\dim(N)$.
In particular, if the codimension $n-m$
is greater than 2, then the connectivity of \thetag{$*\!*\!*$} tends to
infinity with $k$. The suggested interpretation of this result is that,
if $n-m>2$, then
$V\mapsto\emb(V,N)$ behaves more and more like a 
homotopy sheaf on $\cO$,
with respect to $\cJ_k$, as $k$ tends to infinity.

Suppose now that $M\subset N$, so that $\emb(V,N)$ is a based space
for each open $V\subset M$. Then the following general method for calculating 
or partially calculating $\emb(M,N)$ is second to none. 
Try to determine the cofunctors
$$V\mapsto\text{ homotopy fiber of }[T_k\emb(V,N)\to T_{k-1}\emb(V,N)]$$
for the first few $k>0$. These cofunctors admit a surprisingly 
simple description in terms of configuration spaces;
see Theorem 8.5, and \cite{23}. Try to determine
the extensions (this tends to be very hard) and finally specialize, 
letting $V=M$. This program is already outlined in Goodwillie's expanded 
thesis 
\cite{9, \S Intro.C} 
for spaces of concordance embeddings (a special case of a relative case),
with a pessimistic note added in 
revision: ``\dots it might never be [written up] ...''.
It is also carried out to some 
extent in a simple case in \cite{23}. More details on the same case 
can be found in Goodwillie--Weiss \cite{13}.

\remark{{Organization}} Part I (this paper) is about the series 
of cofunctors \linebreak $V\mapsto T_k\emb(V,N)$, called the {\sl Taylor series}
of the cofunctor $V\mapsto\emb(V,N)$. It is also about Taylor series
of other cofunctors of a similar type, but it 
does not address convergence questions. These will 
be the subject of Part II (\cite{13}, joint work with Goodwillie).
\endremark

\remark{{Convention}} Since homotopy limits are so ubiquitous in this paper,
we need a ``convenient'' category of topological spaces with good 
homotopy limits. The category of fibrant simplicial sets is such
a category. In the sequel, ``Space'' with a capital S means
{\sl fibrant simplicial set}. As a rule, we work with (co--)functors whose 
values are Spaces and whose arguments are spaces (say, manifolds). 
However, there are some situations, for example in \S9, where it is a good idea
to replace the manifolds by their singular simplicial sets. Such
a replacement is often understood. 
\endremark
\eject
\section{Good Cofunctors}
\definition{1.1 Definition} A smooth codimension zero embedding 
$i_1\co V\to W$ between smooth manifolds without boundary 
is an {\sl isotopy equivalence}
if there exists a smooth embedding $i_2\co W\to V$ such that
$i_1i_2$ and $i_2i_1$ are smoothly isotopic to $\id_W$ and $\id_V$,
respectively.
\enddefinition

In the sequel $M$ is a smooth manifold without boundary, and $\cO$ is the 
poset of open subsets of $M$, ordered by inclusion. Usually we think of 
$\cO$ as a category, with exactly one morphism $V\to W$ if $V\subset W$,
and no morphism if $V\not\subset W$.   
A cofunctor
(=contravariant functor) $F$ from $\cO$ to Spaces is {\sl good} if it 
satisfies the following conditions.  
\roster
\itemold{(a)} $F$ takes isotopy equivalences to homotopy equivalences.
\itemold{(b)} For any sequence $\{V_i\mid i\ge0\}$ of objects in $\cO$ with
$V_i\subset V_{i+1}$ for all $i\ge0$, the following canonical map
is a weak homotopy equivalence:
$$F(\cup_iV_i)\la\holim_i F(V_i)\,.$$
\endroster

\definition{1.2 Notation} $\cF$ is the category of all good cofunctors from
$\cO$ to Spaces. The morphisms in $\cF$  are the natural transformations.
A morphism 
$g\co F_1\to F_2$ is an {\sl equivalence} if $g_V\co F_1(V)\la F_2(V)$ is a
homotopy equivalence for all $V$ in $\cO$. Two objects in $\cF$ are 
{\sl equivalent} if they can be related by a chain of equivalences. 
\enddefinition

\definition{1.3 Examples} For any smooth manifold $N$ without boundary,
there are cofunctors from $\cO$ to Spaces
given by $V\mapsto\emb(V,N)$ (Space of smooth embeddings)
and $V\mapsto\imm(V,N)$ (Space of smooth immersions). To be more precise,
we think of $\emb(V,N)$ and $\imm(V,N)$ as geometric realizations of 
simplicial sets:  for example, a 0--simplex of $\imm(V,N)$ is a 
smooth immersion $V\to N$, and a 1--simplex in $\imm(V,N)$ is a smooth 
immersion $V\shtimes\Delta^1\to N\shtimes\Delta^1$ respecting the 
projection to $\Delta^1$. 
\enddefinition

\proclaim{1.4 Proposition} {\sl The cofunctors
$\imm(\bla,N)$ and $\emb(\bla,N)$ are good.}
\endproclaim

Part (a) of goodness is easily verified
for both $\imm(\bla,N)$ and $\emb(\bla,N)$. Namely, suppose that
$i_1\co V\to W$ is an 
isotopy equivalence between smooth manifolds, with isotopy inverse
$i_2\co W\to V$ and isotopies $\{h_t\co V\to V\}$, $\{k_t\co W$ $\to W\}$ from
$i_2i_1$ to $\id_V$ and from $i_1i_2$ to $\id_W$, respectively. Then
$\{h_t\co V\to V\}$ gives rise to a map of simplicial sets
$$\imm(V,N)\shtimes\Delta^1\to\imm(V,N)$$
which is a homotopy from $(i_2i_1)^*$ to the identity. Similarly
$\{k_t\co W\to W\}$ gives rise to a homotopy connecting $(i_1i_2)^*$ and the
identity on $\imm(W,N)$. Therefore $\imm(\bla,N)$ is isotopy invariant.
The same reasoning applies to $\emb(\bla,N)$.

To establish part (b) of goodness, we note that it is enough to
consider the case where $M$ is connected. Then a sequence $\{V_i\}$ as 
in part (b) will either be stationary, in which case we 
are done, or almost all the $V_i$ are open manifolds (no compact components).

\proclaim{1.5 Lemma} {\sl Suppose that $V\in\cO$ has no compact components.
Suppose also that $V=\cup_i K_i$ where each $K_i$ is a smooth compact 
manifold with boundary, contained in the interior of $K_{i+1}$, for $i\ge0$.
Then the canonical maps
$$\imm(V,N)\to\holim_i\imm(K_i,N)\,,\qquad
\emb(V,N)\to\holim_i\emb(K_i,N)$$
are homotopy equivalences.}
\endproclaim

\demo{{\bf Proof}} By the isotopy extension theorem, the restriction
from $\emb(K_{i+1},N)$ to $\emb(K_i,N)$ is a Kan fibration of simplicial 
sets. It is a standard result of immersion theory, much more difficult to
establish than the isotopy extension theorem, that the restriction map
from $\imm(K_{i+1},N)$ to $\imm(K_i,N)$ is a Kan fibration. See
especially Haefliger--Poenaru \cite{15}; although this is written in PL language, it 
is one of the clearest references.

Let $\emb_!(V,N)$ be the Space of
{\sl thick} embeddings $V\to N$, that is, embeddings
$f\co V\to N$ together with a {\sl sober} extension of $f$ to an embedding 
$D(\nu_f)\to N$, where $D(\nu_f)$ is the total space of the normal disk
bundle of $f$. (The word {\sl sober} means that the resulting bundle 
isomorphism between the normal bundle of $M$ in $D(\nu_f)$ and $\nu_f$
itself is the canonical one.) Define $\emb_!(K_i,N)$ 
similarly. In the commutative diagram
$$\CD
\emb_!(V,N) @>=>> \lim_i\emb_!(K_i,N) \\
@VV\text{forget} V   @VV\text{forget} V \\
\emb(V,N) @>>> \lim_i\emb(K_i,N)
\endCD$$
the left--hand vertical arrow is a homotopy equivalence by inspection, 
and the right--hand vertical arrow is a homotopy equivalence because,
according to Bousfield--Kan \cite{1}, the canonical map from 
the limit to the homotopy limit of a tower of Kan fibrations is a homotopy
equivalence of simplicial sets. (Hence the limits in the right--hand column 
could be replaced by homotopy limits.) Hence the lower horizontal arrow is 
a homotopy equivalence. {\sl Note:} the lower horizontal arrow is not
always an isomorphism of simplicial sets---injective immersions are not
always embeddings. \qed
\enddemo

Suppose now that $V=\cup_iV_i$ as in part (b) of goodness, and that
$V$ has no compact components. Each $V_i$
can be written as a union $\cup_jK_{ij}$ where each $K_{ij}$ is smooth
compact with boundary, and $K_{ij}$ is contained in the interior of
$K_{i(j+1)}$. Moreover we can arrange that $K_{ij}$ is also 
contained in the interior of $K_{(i+1)j}$. Writing $F(\bla)$ to mean
$\imm(\bla,N)$, we have a commutative diagram of restriction maps
$$\CD
F(V)@>>>\holim_i F(V_i) \\
@VVV              @VVV \\
\holim_i F(K_{ii})@<<< \holim_i \holim_j F(K_{ij})
\endCD$$
where the vertical arrows are homotopy equivalences by 1.5 and the lower
horizontal arrow is a homotopy equivalence by \cite{4, 9.3}.
(Here we identify $\holim_i\holim_j$ with $\holim_{ij}$.) This shows that
the cofunctor $\imm(\bla,N)$ has property (b).
The same argument applies to the cofunctor $\emb(\bla,N)$. Hence 1.4
is proved. \qed

\section{Polynomial Cofunctors}
The following, up to and including Definiton 2.2, is a quotation from
\cite{23}. Suppose that $F$ belongs to $\cF$ and that $V$ belongs to $\cO$, 
and let $A_0,A_1,\dots,A_k$
be pairwise disjoint closed subsets of $V$. Let 
$\cP_{k+1}$ be the power set of $[k]=\{0,1,\dots,k\}$. This is a poset,
ordered by inclusion. We make a functor
from $\cP_{k+1}$ to Spaces by
$$S\mapsto F\bigl(V\minus\cup_{i\in S}A_i\bigr)\tag{$*$}$$   
for $S$ in $\cP_{k+1}$. Recall that, in general, a functor 
$\cX$ from $\cP_{k+1}$ to Spaces is called a $(k+1)$--cube of Spaces.

\definition{2.1 Definition} (\cite{6}, \cite{7})\qua 
The {\sl total homotopy fiber}
of the cube $\cX$ is the homotopy fiber of the canonical map
$$\cX(\emptyset)\la\holim_{S\ne\emptyset}\,\cX(S)\,.$$
If the canonical map $\cX(\emptyset)\to\holim_{S\ne\emptyset}\cX(S)$
is a homotopy equivalence, then $\cX$ is {\sl homotopy Cartesian} or just 
{\sl Cartesian}. 

A {\sl cofunctor} $\cY$ from $\cP_{k+1}$ to spaces will also be called 
a {\sl cube} of spaces, since $\cP_{k+1}$ is isomorphic
to its own opposite. The {\sl total homotopy fiber} of $\cY$ is the homotopy 
fiber of $\cY([k])\to\holim_{S\ne[k]}\cY(S)$.
\enddefinition

Inspired by \cite{7, 3.1} we decree:

\definition{2.2 Definition} The cofunctor $F$ is {\sl polynomial of degree
$\le k$} if the $(k+1)$--cube \thetag{$*$} is Cartesian for arbitrary
$V$ in $\cO$ and pairwise disjoint closed subsets $A_0,\dots, A_k$ of $V$.
\enddefinition

\remark{{Remark}} In Goodwillie's calculus of functors, a functor
from spaces to spaces is {\sl of degree $\le k$} if it takes 
strongly cocartesian $(k+1)$--cubes to Cartesian $(k+1)$--cubes.
The pairwise disjointness condition in 2.2 is there precisely to ensure
that the cube given by $S\mapsto V\minus\cup_{i\in S}A_i$ is strongly 
cocartesian. 
\endremark

\definition{2.3 Example} The cofunctor $V\mapsto\imm(V,N)$ 
is polynomial of degree $\le1$ if either $\dim(N)>\dim(M)$ or 
the dimensions are equal and $M$ has no compact component. 
This amounts to saying that for open 
subsets $V_1$ and $V_2$ of $M$, the following square of restriction maps
is a homotopy pullback square:
$$\CD
\imm(V_1\cup V_2,N)@>>>\imm(V_1,N) \\
@VVV                   @VVV \\
\imm(V_2,N) @>>> \imm(V_1\cap V_2,N).
\endCD$$
To prove this we use lemma 1.5. Accordingly it is enough to prove that
$$\CD
\imm(K_1\cup K_2,N)@>>>\imm(K_1,N) \\
@VVV                   @VVV \\
\imm(K_2,N) @>>> \imm(K_1\cap K_2,N)   
\endCD\tag{$**$}$$
is a homotopy pullback square whenever $K_1,K_2\subset M$ are smooth
compact codimension zero submanifolds of $M$ whose boundaries intersect
transversely. (Then $K_1\cap K_2$ is smooth ''with corners''.) But 
\thetag{$**$} is a strict pullback square of Spaces
in which all arrows are 
(Kan) fibrations, by \cite{15}. \qed
\enddefinition

\definition{2.4 Example} Fix a space $X$, and for $V\in\cO$ let 
$\binom Vk$ be the configuration space of unordered $k$--tuples in $V$.
This is the complement of the fat diagonal in the $k$--fold symmetric product
$(V\shtimes V\shtimes\dots\shtimes V)/\Sigma_k$. The cofunctor
$$V\mapsto\map\left(\binom Vk,X\right)$$
where $\map$ denotes a simplicial set of maps, is polynomial of degree $\le k$.
Here is a sketch proof:  Let $A_0$, $A_1$, \dots, $A_k$ be pairwise 
disjoint closed subsets of $V$. Any unordered $k$--tuple in $V$ must have
empty intersection with one of the $A_i$. Therefore
$$\binom Vk=\bigcup_i\binom{V\minus A_i}k$$
and it is not hard to deduce that the canonical map
$$\hocolim\Sb S\subset\{0,1,\dots,k\} \\ S\neq\emptyset\endSb
\binom{V\minus\cup_{i\in S}A_i}k\quad\la\quad \binom Vk$$
is a homotopy equivalence. Compare \cite{24}. Applying $\map(\bla,X)$ turns the 
homotopy colimit into a homotopy limit and the proof is complete. \qed 
\enddefinition

\definition{2.5 Example} Let $\cA$ be a small category and let
$\phi\co \cA\to\cF$ be a functor, which we will write in the form
$a\mapsto\phi_a$. Suppose that each $\phi_a$ is polynomial of degree
$\le k$. Then 
$$V\mapsto \holim_a\phi_a(V)$$
is in $\cF$, and is polynomial of degree $\le k$. Special case:
For $\cA$ take the poset of nonempty subsets of $\{0,1\}$, and conclude
that $\cF$ is closed under homotopy pullbacks.
\enddefinition
 
\section{Special Open Sets}
Let $\cO k$ consist of all open subsets of $M$ which are 
diffeomorphic to a disjoint union of at most $k$ copies of $\RR^m$,
where $m=\dim(M)$. We think of $\cO k$ as a full subcategory of $\cO$.
There is an important relationship between $\cO k$ and definition
2.2 which we will work out later, and which is roughly as follows.
A good cofunctor from $\cO$ to Spaces
which is polynomial of degree $\le k$ is {\sl determined} by its
restriction to $\cO k$, and moreover the restriction to $\cO k$ can be 
{\sl arbitrarily prescribed}.--- In this section, 
however, we merely examine the homotopy type of $|\cO k|$
and use the results to study the process of {\sl inflation}
(right Kan extension) of a cofunctor
along the inclusion $\cO k\hookrightarrow\cO$. 

For the proof of lemma 3.9 below, we need {\sl double categories}
\cite{17}. Recall first
that a category $\cC$ consists of two classes, $ob(\cC)$ and 
$mor(\cC)$, as well as maps $s,t\co mor(\cC)\to ob(\cC)$ ({\sl source} and 
{\sl target}) and $1\co ob(\cC)\to mor(\cC)$ and
$$\circ\co mor(\cC)_t\shtimes_smor(\cC)\la mor(\cC)$$
({\sl composition}), where $_t\shtimes_s$ denotes the fibered product
(or pullback) over $ob(\cC)$. The maps $s,t,1$ and $\circ$ satisfy certain
relations.
A {\sl double category} is a category object in the 
category of categories. Thus a double category $\cC$ consists of two categories,
$ob(\cC)$ and $mor(\cC)$, as well as functors
$s,t\co mor(\cC)\to ob(\cC)$ ({\sl source} and 
{\sl target}) and $1\co ob(\cC)\to mor(\cC)$ and 
$$\circ\co mor(\cC)_t\shtimes_smor(\cC)\la mor(\cC)$$
({\sl composition}) where $_t\shtimes_s$ denotes the fibered product
(or pullback) over $ob(\cC)$. These functors $s,t,1$ and $\circ$ satisfy the 
expected relations. Alternative definition:  A double category consists of 
four classes, $ob(ob(\cC))$, $mor(ob(\cC))$, $ob(mor(\cC))$ and 
$mor(mor(\cC))$, and certain maps relating them \dots This definition
has the advantage of being more symmetric. In particular, we see that
a double category $\cC$ determines two ordinary categories, 
the {\sl horizontal} category $\cC_h$ and the {\sl vertical} category
$\cC_v$, both with object class $ob(ob(\cC))$. The morphism class of
$\cC_h$ is $ob(mor(\cC))$, that of $\cC_v$ is $mor(ob(\cC))$.

The {\sl nerve} of a double category $\cC$ is a 
bisimplicial set, denoted by $|\cC|$.
  
\definition{3.1 Example} Suppose that two groups $H$ and $V$ 
act on the same set $S$ (both on the left). Make a double category
$\cC$ with $ob(ob(\cC))=S$, $ob(mor(\cC))\!=\!S\shtimes H$,
$mor(ob(\cC))=S\shtimes V$, and
$$mor(mor(\cC)) :=\{(s,h_1,h_2,v_1,v_2)\mid v_2h_1s=h_2v_1s\}\,.$$
Thus an element in $mor(mor(\cC))$ is a ''commutative diagram''
$$\CD
v_1s@>h_2>>h_2v_1s=v_2h_1s \\
@AA v_1 A    @AA v_2 A \\
s @>h_1>> h_1s
\endCD$$
where the vertices are in $S$ and the labelled arrows indicate left
multiplication by suitable elements of $H$ or $V$. 
\enddefinition

\definition{3.2 Example} An ordinary category $\cA$ gives rise to
a double category denoted $\cA\cA$ with $(\cA\cA)_h=\cA=(\cA\cA)_v$ and 
with $mor(mor(\cA\cA))$ equal to the class of commutative 
squares in $\cA$. More generally, if $\cA$ is a subcategory of 
another category $\cB$ containing all objects of $\cB$, then we can form a double category
$\cA\cB$ such that $(\cA\cB)_h=\cB$, $(\cA\cB)_v=\cA$, and
such that $mor(mor(\cA\cB))$ is the class of commutative squares
in $\cB$ whose vertical arrows belong to the subcategory $\cA$:
$$\CD
C@>>>D \\
@AAA     @AAA \\
A @>>> B \,. 
\endCD$$
\enddefinition
\proclaim{3.3 Lemma {\rm \cite{22, Lemma 1.6.5}}}
{\sl The inclusion of nerves $|\cB|\to|\cA\cB|$ is
a homotopy equivalence.}
\endproclaim 

Recall that the homotopy limit of a cofunctor $F$ 
from a small (ordinary) category $\cC$ to $\cT$, the category of Spaces,
is the totalization of the cosimplicial Space
$$p\mapsto\prod_{G\co [p]\to\cC}F(G(0))$$
where the product is taken over all functors $G$ from 
$[p]=\{0,1,\dots,p\}$ to $\cC$. What can we do if $\cC$ is a 
double category and $F$ is a (double) cofunctor from $\cC$ to $\cT\cT$? Then we 
define the homotopy limit as the totalization of the bi--cosimplicial Space
$$(p,q)\mapsto\prod_{H\co [p]\shtimes[q]\to\cC}F(G(0,0))\,.$$
Note that $[p]\shtimes[q]$ is a double category, horizontal arrows being
those which do not change the second coordinate and vertical arrows being those
which do not change the first coordinate. 

We need a variation on 3.3 involving homotopy limits. In the situation of 3.3,
assume that $F$ is a cofunctor from $\cB$ to Spaces (=$\cT$)
taking all morphisms 
in $\cA$ to homotopy equivalences. We can think of $F$ as a 
double cofunctor from
$\cA\cB$ to $\cT\cT$.

\proclaim{3.4 Lemma} {\sl The projection
$$\holim_{\cA\cB}F\to \holim_{\cB}F$$
is a homotopy equivalence.}
\endproclaim

\demo{{\bf Proof}} Let $\cA_p\cB$ be the ordinary category whose objects are diagrams
of the form $A_0\to\dots\to A_p$
in $\cA$, with natural transformations in $\cB$ between such diagrams as 
morphisms. It is enough to show that the face functor
$$d\co (A_0\to\dots\to A_p)\mapsto A_0$$
induces a homotopy equivalence
$$d^*\co \holim_{\cB} F\la \holim_{\cA_p\cB} Fd\,.$$ 
The face functor $d$ has an obvious left adjoint, say $e$. Thus there is 
a natural transformation $\tau$ from $ed$ to the identity on $\cA_p\cB$.
The natural transformation is a functor
$$\tau\co [1]\shtimes\cA_p\cB\la\cA_p\cB\,.$$
Now the key observation is that $Fd\tau$ equals the composition
$$[1]\shtimes\cA_p\cB@>\text{ projection }>>\cA_p\cB@>Fd>>\cT\,.$$
Hence $\tau^*$ can be defined as a map from
$\holim Fd$ to $\holim (Fd\cdot\text{proj})$. Now 
$i_0^*\tau^*=(ed)^*$ and $i_1^*\tau^*=\id$, where $i_0$ and $i_1$
are the standard injections of $\cA_p\cB$ in $[1]\shtimes\cA_p\cB$.
Therefore $(ed)^*$ is homotopic to the identity. Also, $de$ is an
identity functor. \qed
\enddemo

To be more specific now, let $\cI k\subset\cO k$ be the subcategory consisting
of all morphisms which are isotopy equivalences. Eventually we will be interested in the 
double category $\cI k\cO k$. Right now we need a lemma concerning $\cI k$ itself.
\proclaim{3.5 Lemma} 
$$|\cI k|\simeq\coprod\limits_{0\le j\le k}\binom Mj\,.$$
\endproclaim

\demo{{\bf Proof}} Observe that $\cI k$ is a coproduct 
$\coprod\cI^{(j)}$ where $0\le j\le k$ and the objects of $\cI^{(j)}$
are the open subsets of $M$ diffeomorphic to a union of $j$ copies of
$\RR^m$. We have to show
$$|\cI^{(j)}|\simeq\binom Mj\,.$$
For $j=0$ this is obvious. Here is a proof for $j=1$, following \cite{5, 3.1}. 
Let $E\subset|\cI^{(1)}|\shtimes M$ consist of all pairs
$(x,y)$ such that the (open) cell of $|\cI^{(1)}|$ containing $x$
corresponds to a nondegenerate simplex (diagram in $\cI^{(1)}$)
$$V_0\to V_1\to \dots\to V_r$$
where $y\in V_r$. The projection maps
$$|\cI^{(1)}|@<<<E@>>>M$$ 
are {\it almost locally trivial} in the sense of \cite{20, A.1}, since $E$ 
is open in $|\cI^{(1)}|\shtimes M$. By 
\cite{20, A.2} it is enough to verify that both have contractible
fibers. Each fiber of $E@>>>|\cI^{(1)}|$ is homeomorphic to 
euclidean space $\RR^n$. 

Let $E_y$ be the fiber of $E\to M$ over $y\in M$. 
This embeds in $|\cI^{(1)}|$ under the projection, and we can describe it as the union 
of all open cells corresponding to nondegenerate simplices
$(U_0\to\dots\to U_k)$ where $U_k$ contains $y$. There is a subspace
$D_y\subset E_y$ defined as the union of all
open cells corresponding to nondegenerate simplices
$(U_0\to \dots\to U_k)$ where $U_0$ contains $y$. 
Note the following:
\roster
\itemold{$\bullet$} $D_y$ is a deformation retract of $E_y$. Namely, 
suppose that $x$ in $E_y$ belongs to a cell corresponding to a simplex
$(U_0,\dots,U_k)$
with $y\in U_k$. Let $(x_0,x_1,\dots,x_k)$ be the barycentric
coordinates of $x$ in that simplex, all $x_i>0$, and let $j\le k$ be the 
least integer such
that $y\in U_j$. Define a deformation retraction by 
$$\aligned
h_{1-t}(x) :=&(tx_{\text{no}}+
x_{\text{yes}})^{-1}(tx_0,\dots,tx_{j-1},x_j,\dots,x_k)
\\ &x_{\text{no}} :=\sum_{i<j}x_i\,\qquad x_{\text{yes}} :=
\sum_{i\ge j}x_i\endaligned$$
for $t\in[0,1]$, using the barycentric coordinates in the same simplex.
\itemold{$\bullet$} $D_y$ is homeomorphic to the classifying space of 
the poset of all $U\in\cI^{(1)}$ containing $y$. The opposite poset 
is directed, so $D_y$ is contractible. 
\endroster
Hence $E_y$ is contractible, and the proof for $j=1$ is complete. 
In the general case  $j\ge 1$ let 
$$E \subset \cI^{(j)}\shtimes\binom Mj$$ 
consist of all pairs
$(x,S)$ such that the (open) cell of $|\cI^{(j)}|$ containing $x$
corresponds to a nondegenerate simplex  
$$V_0\to V_1\to \dots\to V_r$$
(diagram in $\cI^{(j)}$) 
where each component of $V_r$ contains exactly one point from $S$.
Again the projections from $E$ to $|\cI^{(j)}|$ and to $\binom Mj$ are 
homotopy equivalences. \qed
\enddemo

For $p\ge0$ let $\cI k\cO k_p$ be the category whose objects are 
functors $G\co [p]\to\cO k$ and whose morphisms are double functors
$$[1]\shtimes[p]\la \cI k\cO k\,.$$
(Note that the nerve of the simplicial category $p\mapsto \cI k\cO k_p$
is isomorphic to the nerve of the double category $\cI k\cO k$.)
The rule $G\mapsto G(p)$ is a functor from $\cI k\cO k_p$ to $\cI k$.
In the next lemma we have to make explicit reference to $M$ and another
manifold $V$, so we write $\cO k(M)$, $\cI k(M)$ and so on. 
\proclaim{3.6 Lemma} {\sl For any object $V$ in $\cI k(M)$, the homotopy 
fiber over the 0--simplex $V$ of the map
$$|\cI k\cO k_p(M)|\la |\cI k(M)|$$
induced by $G\mapsto G(p)$ is homotopy equivalent 
to $|\cI k\cO k_{p-1}(V)|$.}
\endproclaim

\definition{3.7 Remark} Combining 3.6 and 3.5, and induction on $p$,
 we can get a very good idea of the homotopy type of $|\cI k\cO k_p(M)|$.
In particular, the functor
$$V\mapsto |\cI k\cO k_p(V)|$$
from $\cO=\cO(M)$ to Spaces takes isotopy equivalences to homotopy
equivalences because the functors $V\mapsto\binom Vj$ have this property.
\enddefinition

\demo{Proof of 3.6} Using Thomason's homotopy colimit theorem 
\cite{21}
we can make the identification
$$|\cI k\cO k_p(M)|\simeq\hocolim_{V\in\cI k(M)}|\cI k\cO k_{p-1}(V)|\,.$$
Then the map under investigation corresponds to the projection from
the homotopy colimit to the nerve of $\cI k(M)$. This map is already a 
quasifibration of simplicial sets. Namely, all morphisms $V_1\to V_2$
in $\cI k(M)$ are isotopy equivalences by definition, and inductively 
we may assume that the functor $V\mapsto|\cI k\cO k_{p-1}(V)|$ takes 
isotopy equivalences to homotopy equivalences (see remark 3.7). Therefore
the homotopy fiber that we are interested in has the same homotopy type as 
the honest fiber. \qed
\enddemo

Let $E$ be a cofunctor from $\cO k=\cO_k(M)$ to Spaces taking 
morphisms in $\cO k$ which are isotopy equivalences to homotopy equivalences.
Use this to define a cofunctor $E^!$ from $\cO$ to Spaces by the formula
$$E^!(V)=\holim\Sb U\in\cO k(V)\endSb E(U)\,.$$
In categorical patois: $E^!$ is the {homotopy right Kan extension} of $E$ 
along the inclusion functor $\cO k\to\cO$.
 
\proclaim{3.8 Lemma} {\sl $E^!$ is good.}
\endproclaim
 
\demo{{\bf Proof}} From 3.4 we know that the projection
$$\holim\Sb U\in\cI k\cO k(V)\endSb E(U)
\la\holim\Sb U\in\cO k(V)\endSb E(U)$$
is a homotopy equivalence. The domain of this projection can be thought of as 
the totalization of the cosimplicial Space
$$p\mapsto\holim_{U_0\to\dots\to U_p} E(U_0)$$
where the homotopy limit, $\holim E(U_0)$, 
is taken over $\cI k\cO k_p(V)$ as defined just
before 3.6. Note that the cofunctor $(U_0\to\cdots\to U_p)\mapsto E(U_0)$
takes all morphisms to homotopy equivalences. Hence its homotopy colimit 
is quasifibered over the nerve of the indexing category, and its 
homotopy limit may be identified (up to homotopy equivalence) with the 
section Space of the associated fibration. Using 3.6 and 3.7 now we see that
$$V\mapsto\holim_{U_0\to\dots\to U_p} E(U_0)$$ 
is a good cofunctor $E^!_p$ for each $p$. Hence $E^!$ is
good, too. \qed
\enddemo
 
We come to the main result of the section. It is similar to certain well--known statements about
{\sl small simplices}, for example \cite{2, III.7.3},
which are commonly used to prove excision
theorems. Let $\varepsilon$ be an open cover of $M$. We say that 
$V\in\cO k$ is $\varepsilon$--small if each connected component of $V$
is contained in some open set of the cover $\varepsilon$. Let
$\varepsilon\cO k=\varepsilon\cO k(M)$ be the full 
sub--poset of $\cO k$ consisting of
the $\varepsilon$--small objects. For $V\in\cO$ let 
$$\varepsilon E^!(V) :=\holim_{U\in\varepsilon\cO k(V)}E(U)\,.$$

\proclaim{3.9 Theorem} {\sl The projection $E^!(V)\to\varepsilon E^!(V)$ 
is a homotopy equivalence.}
\endproclaim

\demo{{\bf Proof}} Using the notation from the proof of 3.8, and obvious 
$\varepsilon$--modifications, we see that it suffices to prove that
the projection $E^!_p(V)\to \varepsilon E^!_p(V)$ is a homotopy equivalence,
for all $V$ and $p$. However, the analysis of $E^!_p(V)$ as a section Space
(proof of 3.8) works equally well for $\varepsilon E^!_p(V)$, and gives the
same result up to homotopy equivalence. In particular
3.5 and 3.6 go through in the $\varepsilon$--setting. \qed
\enddemo

\section{Construction of Polynomial Cofunctors} 

We continue to assume that $E$ is a cofunctor from $\cO k$ to Spaces 
taking isotopy equivalences to homotopy equivalences. 

\proclaim{4.1 Theorem} {\sl The cofunctor $E^!$ on $\cO$ is polynomial of
degree $\le k$.}
\endproclaim

\demo{{\bf Proof}} We have to verify that the condition in 2.2 is satisfied.
Without loss of generality, $V=M$. Then $A_0,A_1,\dots,A_k$
are pairwise disjoint closed subsets of $M$. Let $M_i=M\minus A_i$ and 
$M_S=\cap_{i\in S}M_i$ for $S\subset\{0,1,\dots,k\}$.
Using 3.9, all we have to show is that the $(k+1)$--cube of Spaces 
$$S\mapsto \varepsilon E^!(M_S)$$
is homotopy Cartesian. Here $\varepsilon$
can be any open cover of $M$, and in the present circumstances we choose 
it so that none of the open sets in $\varepsilon$ meets more than 
one $A_i$. Then 
$$\varepsilon\cO k=\bigcup_i\varepsilon\cO k(M_i)\,.$$
(This is the pigeonhole principle again:  Each component of an object
$U$ in $\varepsilon\cO_k$  meets at most one of the $A_i$, but 
since $U$ has at most $k$ components, $U\cap A_i=\emptyset$ for some $i$.)
With lemma 4.2 below, we conclude that the canonical map
$$\holim_{\varepsilon\cO k}E\la\holim_{S\ne\emptyset}
\holim_{\varepsilon\cO k(M_S)}E$$
is a homotopy equivalence. But this is what we had to show. \qed
\enddemo

In lemma 4.2 just below, an {\sl ideal} in a poset $\cQ$ is a subset $\cR$ of 
$\cQ$ such that for every $b\in\cR$, all $a\in\cQ$ with $a\le b$ belong 
to $\cR$. 

\proclaim{4.2 Lemma} {\sl Suppose that the poset $\cQ$ is a 
union of ideals $\cQ_i$, where $i\in T$. For finite nonempty
$S\subset T$ let $\cQ_S=\cap_{i\in S}\cQ_i$. Let $E$
be a cofunctor from $\cQ$ to Spaces. Then the canonical map
$$\holim_{\cQ} E\la\holim_{S\ne\emptyset}\holim_{\cQ_S}E$$
is a homotopy equivalence.}
\endproclaim

\demo{{\bf Proof}} Let $\int\cQ_S$ be the poset consisting of all pairs
$(S,x)$ where $S\subset T$ is finite, nonempty and where $(S,x)\le(S',y)$ if and
only if $S\subset S'$ and $x\le y$ in $\cQ$. The forgetful map taking
$(S,x)$ to $x$ is a functor $u\co \int\cQ_S\to\cQ$. It is 
right cofinal, ie the {\sl under} category $y\down u$ is
contractible for every $y$ in $\cQ$. Therefore, by the {\sl cofinality
theorem for homotopy inverse limits} \cite{1, ch.~XI,~9.2}, 
\cite{4, 9.3} the 
obvious map
$$\holim_{x\in\cQ}E(x)\la\holim_{(S,x) \in\int\cQ_S}E(x)$$
is a homotopy equivalence. (Note that it has to be right cofinal instead
of the usual left cofinal because we are dealing with a {\sl cofunctor}
$E$.) By inspection, the codomain of this map is homeomorphic to
$$\holim_{S\ne\emptyset}\holim_{\cQ_S}E\,.\qed$$
\enddemo

\remark{{Remark}} Note that the obvious map $E(U)\to E^!(U)$ is a homotopy
equivalence for every $U$ in $\cO k$. This is again an application of
the cofinality theorem for homotopy inverse limits, although a much more
obvious one. In this sense $E^!$ extends $E$.
\endremark

\section{Characterizations of Polynomial Cofunctors}

\proclaim{5.1 Theorem} {\sl Let $\gamma\co F_1\to F_2$ be a morphism in
$\cF$. Suppose that
both $F_1$ and $F_2$ are polynomial of degree $\le k$. If 
$\gamma\co F_1(V)\to F_2(V)$
is a homotopy equivalence for all $V\in\cO k$, then it is a homotopy 
equivalence for all $V\in\cO$.}
\endproclaim

\demo{{\bf Proof}} Suppose that $\gamma\co F_1(V)\to F_2(V)$ is a homotopy 
equivalence for all $V\in\cO k$. Suppose also that $W\in\cO r$, 
where $r>k$. Let 
$A_0,A_1,\dots,A_k$ be distinct components of $W$ and let
$W_S=\cap_{i\in S}(W\minus A_i)$ for $S\subset\{0,1,\dots,k\}$. Then 
$$F_i(W)\,\,\simeq\,\,\holim_{S\ne0} F_i(W_S)$$
for $i=1,2$ and therefore $\gamma\co F_1(W)\to F_2(W)$ is a 
homotopy equivalence provided $\gamma$ from $F_1(W_S)\to F_2(W_S)$ is a 
homotopy equivalence for all nonempty $S\subset\{0,1,\dots,k\}$. But 
$W_S$ for $S\ne\emptyset$ has fewer components than $W$, so by induction
the proviso is correct. This takes care of all $W\in\cup_r\cO r$. 

Next, suppose that $W=\intr(L)$ where $L$ is a smooth compact codimension
zero submanifold of $M$. Choose a handle decomposition for $L$, let
$s$ be the maximum of the indices of the handles, and let $t$ be the number
of handles of index $s$ that occur.  If $s=0$
we have $W\in\cO r$ for some $r$ and this case has been dealt with. 
If $s>0$, let  
$e\co \DD^{m-s}\shtimes\DD^s\to L$ be one of the $s$--handles. We assume that 
$e^{-1}(\partial L)$ is $\partial\DD^{m-s}\shtimes\DD^s$. Since $s>0$ we can
find pairwise disjoint small closed disks $C_0,\dots,C_k$ in $\DD^s$
and we let
$$A_i :=e(\DD^{m-s}\shtimes C_i)\cap W$$
for $0\le i\le k$. Then each $A_i$ is closed in $W$ and $W\minus A_i$ is the 
interior of a smooth handlebody in $M$ which has a handle decomposition
with no handles of index $>s$, and fewer than $t$ handles of index $s$. 
The same is true for $W_S :=\cap_{i\in S}(W\minus A_i)$ provided 
$S\ne\emptyset$. Therefore, by induction,
$$\gamma\co F_1(W_S)\la F_2(W_S)$$
is a homotopy equivalence for $\emptyset\ne S\subset\{0,1,\dots,k\}$ 
and consequently the right--hand vertical arrow in
$$\CD
F_1(W)@>>>\holim\limits_{S\ne\emptyset} F_1(W_S) \\
@VV\gamma V       @VV\gamma V \\
F_2(W)@>>>\holim\limits_{S\ne\emptyset} F_2(W_S) 
\endCD$$
is a homotopy equivalence. But the two horizontal arrows are also
homotopy equivalences, because $F_1$ and $F_2$ are polynomial of degree 
$\le k$. Therefore the left--hand vertical arrow is a homotopy 
equivalence. This takes care of every $W\in\cO$ which is the interior 
of a compact smooth handlebody in $M$.

The general case follows because $F_1$ and $F_2$ are good cofunctors;
see especially property (b) in the definition of goodness, just after 1.1.
\qed
\enddemo

For $F$ in $\cF$ let $T_kF$ be the homotopy right Kan extension of
the restriction of $F$ to $\cO k$. The explicit formula is 
$$T_kF(V) :=\holim_{U\in\cO k(V)}F(U)\,.$$
From \S3 and \S4 we know that $T_kF$ is good and polynomial of degree $\le k$. 
There is an obvious forgetful morphism $\eta_k\co F\to T_kF$. 
Clearly the natural map 
$\eta_k\co F(U)\to T_kF(U)$ is a homotopy equivalence for every $U\in\cO k$. 
Hence, by 5.1, if $F$ is already polynomial of degree $\le k$, 
then $\eta_k$ from $F(V)$ to $T_kF(V)$ is a homotopy equivalence for every 
$V\in\cO$. In this sense an $F$ which is polynomial of degree $\le k$
is determined by its restriction $E$ to $\cO k$. The restriction
does of course take isotopy equivalences in $\cO k$ to homotopy equivalences.
We saw in \S4 that that is essentially the only condition it must satisfy. 

The polynomial objects in $\cF$ can also be characterized in sheaf
theoretic terms. Recall the Grothendieck topologies $\cJ_k$ on $\cO$, from
the introduction.

\proclaim{5.2 Theorem} {\sl A good cofunctor $F$ from $\cO$ to Spaces is 
polynomial of degree $\le k$ if and only if it is a homotopy sheaf with 
respect to the Grothendieck topology $\cJ_k$.}
\endproclaim

\demo{{\bf Proof}} Suppose that $F$ is a homotopy sheaf with respect to $\cJ_k$.
Let $V\in\cO$ and pairwise disjoint closed subsets $A_0,\dots,A_k$ of $V$ 
be given. Let $V_i=V\minus A_i$. Then the inclusions $V_i$ for 
$0\le i\le k$ form a covering of $V$ in the Grothendieck topology $\cJ_k$.
Hence
$$F(V)\la \holim_{R}F(\cap_{i\in R}V_i)$$
is a homotopy equivalence; the homotopy limit is taken over the 
nonempty subsets $R$ of $\{0,\dots,k\}$. This shows that $F$ is 
polynomial of degree $\le k$. 

Conversely, suppose that $F$ is polynomial of degree $\le k$. Let $W\in\cO$
be given and let $\{V_i\to W\mid i\in S\}$ be a covering of $W$ in the
Grothendieck topology $\cJ_k$. Let $E$ be the restriction of $F$ to $\cO k$.
Define $\varepsilon E^!$ as in \S3, just before 3.9, where $\varepsilon$ 
is the covering $\{V_i\}$. Up to equivalence, $F$ and $\varepsilon E^!$ are
the same. By 4.2, the canonical map
$$\varepsilon E^!(W) \la \holim_R\varepsilon E^!(\cap_{i\in R}V_i)$$
is a homotopy equivalence. Here again, $R$ runs through the finite 
nonempty subsets $R$ of $S$. \qed
\enddemo
 
\section{Approximation by Polynomial Cofunctors} 

From \S5, we have for every $k\ge 0$ an 
endofunctor $T_k\co \cF\to\cF$ given by the rule $F\mapsto T_kF$, and a natural
transformation from the identity $\cF\to\cF$ to $T_k$ given by 
$\eta_k\co F\to T_kF$ for all $F$. It is sometimes convenient to
define $T_{-1}$ as well, by $T_{-1}F(V) :=*$.
The following theorem is mostly a summary of results from \S5. It 
tries to say that $T_k$ is essentially {\sl left adjoint} to the
inclusion functor $\cF_k\to\cF$. Here $\cF_k$ is the full subcategory
of $\cF$ consisting of the objects which are polynomial of degree
$\le k$. Compare \cite{25, Thm.6.1}.

\proclaim{6.1 Theorem} {\sl The following holds for every $F$ in $\cF$ and 
every $k\ge0$.
\roster
\item $T_kF$ is polynomial of degree $\le k$.
\item If $F$ is already polynomial of degree $\le k$, then  
$\eta_k\co F\to T_kF$ is an equivalence.
\item $T_k(\eta_k)\co T_kF\to T_k(T_kF)$ is an equivalence.
\endroster}
\endproclaim

\demo{{\bf Proof}} Properties \therosteritem1 and \therosteritem2 have been 
established in \S5. As for \therosteritem3, we can use 5.1 and we then only 
have to verify that
$$T_k(\eta_k)\co T_kF(W)\to T_k(T_kF(W))$$
is a homotopy equivalence for every $W\in\cO k$. Written out in detail the
map takes the form
$$\align \holim_{V\in\cO k(W)}F(V)\la&\holim_{V\in\cO k(W)}T_kF(V) \\
              =&\holim_{V\in\cO k(W)}\quad\holim_{U\in\cO k(V)}F(U)
\endalign$$
and it is induced by the maps $F(V)\to\holim_U F(U)$ for $V$ in 
$\cO k(V)$. These maps are 
clearly homotopy equivalences, since the identity morphism $V\to V$
is a terminal object in $\cO k(V)$. \qed
\enddemo

\remark{{Remark}} One way of saying that the inclusion of a full subcategory,
say $\cA\to\cB$, has a left adjoint is to say that there exists 
a functor $T\co \cB\to\cB$ and a natural transformation $\eta\co \id_{\cB}\to T$
with the following properties.
\roster
\item $T(b)$ belongs to $\cA$ for every $b$ in $\cB$.
\item For $a$ in $\cA$, the morphism $\eta\co a\to T(a)$ is an isomorphism.
\item For $b$ in $\cB$, the morphism $T(\eta)\co T(b)\to T(T(b))$ is an 
isomorphism.
\endroster
\endremark

From the definitions, there are forgetful 
transformations $r_k\co T_kF\to T_{k-1}F$ for any $F$ and any $k>0$.
They satisfy the relations $r_k\eta_k=\eta_{k-1}\co F\to T_{k-1}F$.
Therefore
$$\{\eta_k\}\co F\la\holim_k T_kF\tag{$*$}$$
is defined. The codomain, with its inverse filtration, may be called the
{\sl Taylor tower} of $F$. Usually one wants to know whether \thetag{$*$}
is a homotopy equivalence. More precisely one can ask two questions:
\roster
\itemold{$\bullet$} Does the Taylor tower of $F$ converge?
\itemold{$\bullet$} If it does converge, does it converge to $F$?
\endroster 
Regarding the first question:  although $\holim_k T_kF$ is always defined,
we would not speak of convergence unless the connectivity of 
$r_k\co T_kF(V)\to T_{k-1}F(V)$ tends to infinity with $k$, independently of $V$.

\section{More Examples of Polynomial Cofunctors} 

\definition{7.1 Example} Let $p\co Z\to\binom Mk$ be a fibration. For
$U\subset\binom Mk$ let $\Gamma(p~;U)$ be the Space of partial sections
of $p$ defined over $U$.  
The cofunctor $F$ on $\cO$ defined by $F(V) :=\Gamma(p~;\binom Vk)$ is 
good and, moreover, it is polynomial of degree $\le k$.  This can be 
proved like 2.4.
\enddefinition

Keep the notation of 7.1. Let 
$\blacktriangle_kV$ be the complement of $\binom Vk$ in 
the $k$--fold symmetric power 
$\spr_kV :=(V\shtimes V\shtimes\dots\shtimes V)/\Sigma^k\,$. The homotopy 
colimit in the next lemma is taken over the poset of all neighborhoods $Q$
of $\blacktriangle_kV$ in $\spr_kV$.

\proclaim{7.2 Lemma} {\sl The cofunctor $G$ on $\cO$ given by
$$G(V) :=\hocolim_Q\Gamma(p~;\tbinom Vk\cap Q)$$
is good.} 
\endproclaim

\demo{{\bf Proof}} We concentrate on part (b) of goodness to begin with.
Fix $V$ and choose a smooth triangulation
on the $k$--fold product $(V)^k$, equivariant with respect to the
symmetric group $\Sigma_k$. Then $\spr_kV$ has a preferred PL structure and 
$\blacktriangle_kV$ is a PL subspace, so we can speak of
{\sl regular neighborhoods} of $\blacktriangle_kV$. It is clear that
all regular neighborhoods of $\blacktriangle_kV$ have the same 
homotopy type, and that each neighborhood of $\blacktriangle_kV$ 
contains a regular one. Therefore, if $L$ is a 
regular neighborhood of $\blacktriangle_kV$, then the canonical
inclusion
$$\Gamma(p~;\tbinom Vk\cap\intr(L))\la\hocolim_Q\Gamma(p~;\tbinom Vk\cap Q)$$
is a homotopy equivalence. This observation tends to simplify matters.
Another observation which tends to complicate matters is that for an open
subset $U$ of $V$ and a regular neighborhood $L$ as above, the 
intersection of $L$ with $\spr_kU$ will usually not be a regular
neighborhood of $\blacktriangle_kU$. However, we can establish
goodness as follows. Suppose that
$$V=\cup_iK_i$$
where each $K_i$ is a smooth compact codimension zero submanifold of $V$,
and $K_i\subset\intr(K_{i+1})$.  
As in the proof of 1.4, it is enough to show that the canonical map
$$G(V)\la\holim_iG(\intr(K_i))$$
is a homotopy equivalence. Abbreviate $\intr(K_i)=V_i$.
Choose a regular neighborhood $L$ of 
$\blacktriangle_kV$ in $\spr_kV$ such that $L\cap\spr_k(K_i)$ is a 
regular neighborhood of $\blacktriangle_k(K_i)$ in $\spr_k(K_i)$ for 
each $i$.
Then it is not hard to see that the inclusion
$$\Gamma(p~;\tbinom{V_i}k\cap\intr(L))\la
\hocolim_R\Gamma(p~;\tbinom{V_i}k\cap R)$$
is a homotopy equivalence, for each $i$. Therefore, in the 
commutative diagram
$$\CD
\Gamma(p~;\binom Vk\cap\intr(L))@>>>\holim_i\Gamma(p~;\binom{V_i}k\cap\intr(L)) \\
@VVV                   @VVV \\
\hocolim_Q\Gamma(p~;\binom Vk\cap Q)@>>>\holim_i
            \hocolim_R\Gamma(p~;\binom{V_i}k\cap R)
\endCD$$
the two vertical arrows are homotopy equivalences. The upper horizontal
arrow is also a homotopy equivalence by inspection. Hence the
lower horizontal arrow is a homotopy equivalence. This completes the proof
of part (b) of goodness.

Proof of part (a) of goodness:  Suppose that $W\hookrightarrow V$ in $\cO$
is an isotopy equivalence. Let $\{j_t\co V\to V\}$ be a smooth isotopy of 
embeddings, with $j_0=\id_V$ and $\im(j_1)=W$. Let
$$X :=\hocolim_R\,\Gamma(j^*p\,;\,(\tbinom Vk\shtimes I)\,\cap\,R)$$
where $I=[0,1]$ and $j^*p$ is the pullback of $p$ under the map 
$$\binom Vk\times I\la \binom Vk\qquad;\qquad(S,t)\mapsto j_t(S)$$
and $R$ runs over the neighborhoods of $\blacktriangle_kV\times I$ in 
$\spr_kV\times I$. {\sl Key observation:} Every $R$ contains a neighborhood
of the form $Q\shtimes I$, where $Q\subset\spr_kV$. This implies that the
restriction maps
$$G(W)@<\rho_W<<X@>\rho_V>>G(V)$$
(induced by the bundle maps $j_1^*p@>>>j^*p@<<<j_0^*p$) are homotopy 
equivalences. The restriction map $G(V)\to G(W)$ that we are 
interested in can be written as a composition
$$G(V)@>j^*>>X@>\rho_W>>G(W)$$
where the arrow labelled $j^*$ is right inverse to $\rho_V$. Therefore 
the restriction map $G(V)\to G(W)$ is a homotopy equivalence. \qed 
\enddemo

\proclaim{7.3 Lemma} {\sl The cofunctor $G$ in 7.2. is polynomial of 
degree $\le k$.}
\endproclaim

\demo{{\bf Proof}} Fix  $W\in\cO$ and let $A_0,\dots,A_k$ be closed and 
pairwise disjoint in $W$. Let $W_i :=W\minus A_i$ and choose 
neighborhoods $Q_i$ of $\blacktriangle_kW_i$ in $\spr_kW_i$. Let
$$\gather W_S=\cap_{i\in S}W_i \\
Q_S=\cap_{i\in S}Q_i\endgather$$
for nonempty $S\subset\{0,1,\dots,k\}$, and $W_{\emptyset}=W$,
$Q_{\emptyset}=\cup_iQ_i$. Then
$$\binom Wk\cap Q_{\emptyset}\quad=\quad\bigcup_i\binom{W_i}k\cap Q_i
\quad\simeq\quad\hocolim_{S\ne\emptyset}\,\,\binom{W_S}k\cap Q_S$$
which shows, much as in the proof of 2.4, that the obvious map
$$\Gamma\left(p~;\tbinom Wk\cap Q_{\emptyset}\right)
\la \holim_{S\ne\emptyset}\Gamma\left(p~;\tbinom{W_S}k\cap Q_S\right)$$
is a homotopy equivalence. We can now complete the proof with
two observations. Firstly, the neighborhoods of $\blacktriangle_kW_S$
of the form $Q_S$, as above, form an {\sl initial} subset \cite{17}
in the poset of all neighborhoods. Secondly, there are situations in
which homotopy inverse limits commute (up to homotopy equivalence) with
homotopy direct limits, and this is one of them. Here we are interested
in a double homotopy limit/colimit of the form
$$\holim_{S\ne\emptyset}\hocolim_{Q_0,\dots,Q_k}\,\,(\bla)$$
where the blank indicates an expression depending on $S$ 
and the $Q_i$ (actually only on the $Q_i$ for $i\in S$). 
Clearly sublemma 7.4 below applies. \qed
\enddemo

\proclaim{7.4 Sublemma} {\sl Let $X$ be a functor from a product 
$\cA\shtimes\cB$ to Spaces, where $\cA$ and $\cB$ are posets. 
Suppose that $\cA$ is finite and that $\cB$ is directed. Then
$$\hocolim_{b\in\cB}\holim_{a\in\cA} X(a,b)\quad\simeq\quad
\holim_{a\in\cA}\hocolim_{b\in\cB}X(a,b)\,.$$}\rm

{{\bf Proof}}\qua Since $\cB$ is a directed poset, 
the homotopy colimits may be replaced
by honest colimits \cite{1}. The universal property of colimits
yields a map
$$\colim_{b\in\cB}\holim_{a\in\cA} X(a,b)
\quad\simeq\quad\holim_{a\in\cA}\colim_{b\in\cB}X(a,b)$$
which is an isomorphism, by inspection. \qed
\enddemo
  
\proclaim{7.5 Proposition} {\sl The cofunctor $G$ in 7.2 and 7.3 is in fact
polynomial of degree $\le k-1$.}
\endproclaim

\demo{{\bf Proof}} We must show that $\eta_k\co G\to T_{k-1}G$ is an 
equivalence. Since $G$ and $T_{k-1}G$ are both polynomial of 
degree $\le k$, it is enough to check that 
$$\eta_k\co G(V)\la T_{k-1}G(V)$$
is an equivalence for every $V\in\cO k$. See 5.1. If $V$ belongs to $\cO r$
for some $r<k$, this is obvious. So we may assume that $V$ has 
exactly $k$ connected components, each diffeomorphic to $\RR^m$. Denote these
components by $A_0,\dots,A_{k-1}$. If we can show that the upper 
horizontal arrow in
$$
\CD
G(V)@>>>\holim\limits_{S\ne\emptyset} G(\cup_{i\notin S}A_i) \\
@VVV                       @VVV \\
T_{k-1}G(V)@>>>\holim\limits_{S\ne\emptyset} T_{k-1}G(\cup_{i\notin S}A_i)
\endCD$$
is a homotopy equivalence, then we are done because the lower horizontal
and the right--hand vertical arrows are homotopy equivalences.
However, this follows in the usual manner (compare proof of 2.4 and
of 7.3) from the observation that
$$\binom Vk\cap Q\,=\,\bigcup_i\binom{V\minus A_i}k\cap Q$$
for sufficiently small neighborhoods $Q$ of $\blacktriangle_kV$
in $\spr_kV$. Notice that the observation as such is new because
this time the closed subsets $A_i$ are $k$ in number, not $k+1$. \qed
\enddemo

We are now in a position to understand the relationship between 
$F$ in 7.1 and $G$ in 7.2. There is an obvious inclusion
$e\co F(V)\to G(V)$, natural in $V$.  
\proclaim{7.6 Proposition} {\sl The morphism $T_{k-1}e\co  T_{k-1}F\to T_{k-1}G$ is an 
equivalence.}
\endproclaim

\demo{{\bf Proof}} By 5.1, it suffices to show that $e\co  F(V)\to G(V)$ is 
a homotopy equivalence for any $V$ which is diffeomorphic to a 
disjoint union of $\ell$ copies of $\RR^m$, where $\ell<k$. 
For such a $V$ choose open subsets
$$V=V_0\supset V_1\supset V_2\supset V_3\supset\dots$$
such that the inclusions $V_{i+1}\to V_i$ are isotopy equivalences,
such that the closure of $V_{i+1}$ in $V_i$ is compact, and such that
$\cap_iV_i$ is a discrete set consisting (necessarily) of $\ell$ points,
one in each component of $V$. In the commutative square
$$\CD
F(V)@>\subset>>\hocolim_i F(V_i) \\
@VV\subset V        @VV\subset V \\
G(V)@>\subset>>\hocolim_i G(V_i)
\endCD\tag{$*$}
$$
the horizontal arrows are now homotopy equivalences because $F$ and $G$
take isotopy equivalences to homotopy equivalences. On the other hand,
suppose that $Q$ is a neighborhood of $\blacktriangle_kV_i$ in
$\spr_kV_i$ for some $i$. Then clearly there exists an integer $j>i$
such that all of $\spr_kV_j$ is contained in $Q$. It follows that
the inclusion of $\hocolim_iF(V_i)$ in 
$$\hocolim_i G(V_i)=
\hocolim_i\hocolim_Q \Gamma\left(p~;\tbinom{V_i}k\cap Q\right)$$
is a homotopy equivalence. Hence all arrows in \thetag{$*$} are homotopy
equivalences. \qed
\enddemo

\section{Homogeneous Cofunctors}
\definition{8.1 Definition} A cofunctor $E$ in $\cF$ is {\sl homogeneous of degree} $k$,
where $k\ge0$, if it is polynomial of degree $\le k$ and if $T_{k-1}E(V)$ 
is contractible for each $V\in\cO$.  
\enddefinition

\remark{{Remark}} The cofunctor given by $E(V)=*$ for all $V$ is homogeneous
of degree $k$ for any $k\ge0$. Conversely, if $E$ is homogeneous of 
degree $k$ and homogeneous of degree $\ell$, where $\ell<k$, then clearly
$E(V)\simeq T_{k-1}E(V)\simeq*\,$.
\endremark

\definition{8.2 Example} Let $F$ in $\cF$ be arbitrary,
and select a point
in $F(M)$, if one exists. Then $T_kF(V)$ is pointed 
for all $V$ and $k$. Therefore a new cofunctor $L_kF$ can be defined by
$$L_kF(V) :=\hofiber\,[T_kF(V)\la T_{k-1}F(V)]\,.$$
It follows from 6.1 that $E$ is homogeneous of degree $k$.
\enddefinition

\definition{8.3 Example} Starting with a fibration $p\co Z\to\binom Mk$,
define $F$ as in 7.1 and define $G$ as in 7.2.
Select a point in $G(M)$. Then
$$E(V) :=\hofiber\,[F(V)@>\subset>>G(V)]$$
is defined. It follows from 7.6 that $E$ is homogeneous of degree $k$.
\enddefinition

Example 8.3 deserves to be studied more. Ultimately $E$ has been 
constructed in terms of the fibration $p$, and a partial section of $p$
defined near the fat diagonal $\blacktriangle_kM$. Is it possible to 
recover $p$ from $E$? In particular, for $S\in\binom Mk$, can we describe the
fiber $p^{-1}(S)$ in terms of $E$? 

Note that $S$ is a subset of $M$
with $k$ elements. Let $V$ be a tubular neighborhood of $S\subset M$,
so that $V$ is diffeomorphic to a disjoint union of $k$ copies of $\RR^m$.
Then $S$ belongs to $\binom Vk\subset\binom Mk$
and therefore we have maps
$$E(V)\la F(V)=\Gamma(p~;\tbinom Vk)@>\text{ evaluation }>>p^{-1}(S)\,.$$

\proclaim{8.4 Proposition} {\sl The composite map $E(V)\to p^{-1}(S)$ is a 
homotopy equivalence.}
\endproclaim

Hence we can indeed describe $p^{-1}(S)$ in terms of $E$, up to 
homotopy equivalence:  namely,
as $E(V)$ for a tubular neighborhood $V$ os $S$ in $M$.

\demo{Proof of 8.4} 
Much as in the proof of 7.6 we choose a sequence of open
subsets
$$V=V_0\supset V_1\supset V_2\supset V_3\supset\dots$$
such that the inclusions $V_{i+1}\to V_i$ are isotopy equivalences,
such that the closure of $V_{i+1}$ in $V_i$ is compact, and such that
$\cap_iV_i=S$. We note that
$$F(V)=\prod_j\Gamma(p~;U_j)$$
where the $U_j$ are the connected components of $\binom Vk$. Among these
components we single out $U_0$, the component containing $S$. It is the only
component whose closure in $\spr_kV$ does not meet $\blacktriangle_kV$.
For the remaining components we can use an idea as in the proof of 7.6,
and find
$$G(V)\quad\simeq\quad\prod_{j\ne0}\Gamma(p~;U_j)\,.$$
Therefore $F(V)\simeq E(V)\shtimes G(V)$ and the composition
$$E(V)\to F(V)\to\Gamma(p~;U_0)\to p^{-1}(S)$$
is a homotopy equivalence. \qed
\enddemo

\remark{{Digression}} Knowing all the fibers of a fibration is not the same as 
knowing the fibration. However, in the present case we can also
``describe'' the entire fibration $p$ in 8.3 in terms of the
cofunctor $E$. Recall from the proof of 3.5 the poset $\cI^{(k)}$.
Its elements are the open subsets of $M$ which are diffeomorphic to
a disjoint union of $k$ copies of $\RR^m$, and for $V,W\in\cI^{(k)}$
we decree $V\le W$ if and only if $V\subset W$ and the inclusion is
an isotopy equivalence. We saw that 
$$|\cI^{(k)}|\simeq\binom Mk\,.$$
Since $\cI^{(k)}\subset\cO$, we can restrict $E$ to $\cI^{(k)}$. The 
restricted cofunctor takes all morphisms to homotopy equivalences, so that
the projection
$$\hocolim_{\cI^{(k)}}E\la|\cI^{(k)}|$$
is a quasifibration. The associated fibration is the one we are looking for.
This motivates the following classification theorem for 
homogeneous cofunctors.
\endremark

\proclaim{8.5 Theorem} {\sl Up to equivalence, all objects in $\cF$ which are 
homogeneous of degree $k$ are of the type discussed in 8.3.}
\endproclaim

\remark{{Outline of proof}}
Of course, the digression just above already gives the idea of the proof,
but we have to proceed a little more cautiously. The plan is:  Given 
$E$, homogeneous of degree but not necessarily defined in terms of some 
fibration, construct the appropriate $F$, polynomial of degree $\le k$,
and a morphism $E\to F$. Then show that $F$ is equivalent to a cofunctor
of type $V\mapsto\Gamma(p~;\tbinom Vk)$.
as in 7.1. This step requires a lemma, 8.6 below.
Finally identify $E$ with the homotopy fiber of the canonical 
morphism from $F$ to $T_{k-1}F$. 
\endremark

\proclaim{8.6 Lemma \cite{3, 3.12}} {\sl Suppose that $Y$ is a functor from a small 
category $\cA$ to the category of Spaces. If $Y$ takes all morphisms 
in $\cA$ to homotopy equivalences, then the canonical projection
$\hocolim_{\cA}Y\to|\cA|$
is a quasifibration. The section Space of the associated fibration
is homotopy equivalent to $\holim_{\cA}X$.}
\endproclaim

\demo{Sketch proof of 8.6} The quasifibration statement is obvious. We 
denote the total Space of the associated fibration by $T$, so that
$\hocolim_{\cA}Y\subset T$
by a homotopy equivalence. For the
statement about the section Space, recall that $\holim Y$ can 
be defined as the Space of natural transformations $\tilde*_{\cA}\to Y$,
where $*_{\cA}$ is the constant functor $a\mapsto*$ on $\cA$, and 
$\tilde*_{\cA}$ is a {\sl CW--functor} weakly equivalent to it (some
explanations below). The standard choice is
$$\tilde*_{\cA}(a) :=|\cA\down a|\,.$$
{\sl CW--functor} refers to a functor with a {\sl CW--decomposition} where the
cells are of the form $\RR^i\shtimes\mor(b,\bla)$ for some $b\in\cA$ and
some $i$. {\sl Weakly equivalent to $*_{\cA}$}  means here that there is an 
augmentation
$\tilde*_{\cA}(a)\to *_{\cA}(a)$, natural in $a$, which is a homotopy 
equivalence for each $a$. In other words, $\tilde*_{\cA}(a)$ is always
contractible.--- 
Suppose now that $X$ is {\sl any} CW--functor from $\cA$ to spaces. There 
are obvious embeddings
$$\nat(X,Y)@>\subset>>\map_{|\cA|}(\hocolim X\,,\,\hocolim Y)
@>\subset>>\map_{|\cA|}(\hocolim X\,,\,T)$$
where $\map_{|\cA|}$ is for Spaces of maps over $|\cA|$.
One shows by induction over the skeletons of $X$ that the 
composite embedding 
is a homotopy equivalence. In particular, this holds  
for $X=\tilde*_{\cA}$. \qed
\enddemo

\demo{Proof of 8.5} Suppose that $E$ in $\cF$ is homogeneous of
degree $k$. Define a cofunctor $F_0$ from $\cO$ to Spaces by
$$F_0(V) :=\holim\limits_{U\in\cI^{(k)}(V)}E(U)\,.$$
Here $\cI^{(k)}(V)\subset\cI^{(k)}$ is the full sub--poset 
consisting of all $U\in\cI^{(k)}$ which are contained in $V$. For the 
meaning of $\cI^{(k)}$, see the digression preceding 8.5. By 8.6,
the cofunctor $F_0$ is equivalent to another cofunctor $F_1$ given 
by a formula of type 
$$F_1(V)=\Gamma(q_V)$$
where $q_V$ is a certain fibration on $|\cI^{(k)}(V)|$. The fibration 
$q_V$ is natural in $W$, in the sense that a morphism
$V\subset W$ in $\cO$ induces a map from the total Space of 
$q_V$ to that of $q_W$, covering the inclusion
$$|\cI^{(k)}(V)|\hookrightarrow|\cI^{(k)}(W)|\,.$$
By inspection, this map of total Spaces maps each fiber of $q_V$ to
the corresponding fiber of $q_W$ by a homotopy equivalence. Hence
$F_1$ is equivalent to the cofunctor $F_2$ given by
$$F_2(V) :=\Gamma\left(q_M;|\cI^{(k)}(V)|\right)\,.$$
Finally we know from 3.5 (and proof) that 
$|\cI^{(k)}(V)|\simeq\binom Vk$, and this can be understood as a chain
of natural homotopy equivalences (natural in $V\in\cO$). It follows
easily that $F_2$ is equivalent to a cofunctor $F_3$ given by
a formula of type
$$F_3(V) :=\Gamma\left(p~;\tbinom Vk\right)$$
where $p$ is a fibration on $\binom Mk$. This is exactly the
kind of cofunctor introduced in \S7, so we now write $F :=F_3$. 
From the definition, $F$ belongs to $\cF$. 
Replacing $E$ by an equivalent cofunctor if necessary, we can assume that
$E$ maps directly to $F$ instead of $F_0$. If $S\in\binom Mk$ and
$V$ is a tubular neighborhood of $S\subset M$, then the 
composition
$$E(V)\la F(V)=\Gamma\left(p~;\tbinom Vk\right)@>\text{ eval. }>>
p^{-1}(S)$$
is a homotopy equivalence, by construction and inspection. This is 
of course reminiscent of 8.4. Now form the commutative square
$$\CD
E@>>>F \\
@VV\eta_{k-1}V  @VV\eta_{k-1}V \\
T_{k-1}E@>>> T_{k-1}F
\endCD\tag{$*$}$$
and recall that $T_{k-1}E(V)$ is contractible for all $V\in\cO$.
Given our analysis of $T_{k-1}F$ in \S7, we can complete the proof
of 8.5 by showing that \thetag{$*$} is homotopy Cartesian.
By 2.5 and 5.1, it suffices to check that
$$\CD
E(V)@>>>F(V) \\
@VV\eta_{k-1}V  @VV\eta_{k-1}V \\
T_{k-1}E(V)@>>> T_{k-1}F(V)
\endCD\tag{$**$}$$ 
is homotopy Cartesian for all $V\in\cO k$. If it happens that 
$V\in\cO r\subset\cO k$ for some $r<k$, then we have $E(V)\simeq*$
by homogeneity and $F(V)\to T_{k-1}F(V)$ is a homotopy equivalence,
by \S5 and \S6. If not, then $V$ has $k$ connected components 
and is a tubular neighborhood of some $S\subset M$, where $S\in\binom Mk$.
Using 8.4 now (and 7.6), and our observation above which
seemed so reminiscent of 8.4, we find that \thetag{$**$} is again
homotopy Cartesian. \qed
\enddemo

\section{The Homogeneous Layers of a Good Cofunctor} 

In this section we work with a fixed $F$ in $\cF$ and a distinguished
element $*\in F(M)$, which we call the base point. Since $M$ is the terminal
object in $\cO$, we may then regard $F$ as a cofunctor from $\cO$ to pointed
Spaces. Define $L_kF$ as in 8.2, and call it {\sl the $k$--th
homogeneous layer of $F$}. According to 8.5, the homogeneous cofunctor $L_kF$
can be classified by some fibration $p\co Z\to\binom Mk$, and a 
partial section of it defined near the fat diagonal $\blacktriangle_kM$.
What does $p$ look like? The answer is implicit in the last section.
Recall that $$\binom Mk\simeq|\cI^{(k)},$$ in the notation of 3.5 and 
sequel. For any $V\in\cI^{(k)}$ with components $V_s$, where $s\in\pi_0(V)$, 
the rule taking a subset $R$ of $\pi_0(V)$ to
the Space $F(\cup_{s\in R}V_s)$ is a $k$--cube of Spaces:
$$R\mapsto F(\cup_{s\in R}V_s)\qquad\qquad(R\subset\pi_0(V))\,.\tag{$*$}$$
As such it has a total homotopy fiber (see 2.1)
which we denote by $\Phi(V)$.
Note that $V\mapsto\Phi(V)$ is a cofunctor from $\cI^{(k)}$
to Spaces taking all morphisms to homotopy equivalences.

\proclaim{9.1 Proposition} {\sl The fibration which classifies $L_kF$ is
the one associated with the quasifibration
$$\hocolim_{V\in\cI^{(k)}}\Phi(V)\la |\cI^{(k)}|\,.$$}
\endproclaim

\remark{{Remark}} Our classifying fibrations on $\binom Mk$ should always 
come with partial sections defined near the fat diagonal.
Note that $\Phi$ is a cofunctor from $\cI^{(k)}$ to 
{\sl pointed} Spaces, so that the (quasi)--fibration in 9.1 does in fact
have a preferred global section. 
\endremark

\demo{Proof of 9.1} Write $j=k-1$ (for typographic reasons).
By \S8, it is enough to show that $L_kF(V)\simeq\Phi(V)$ for 
$V\in\cI^{(k)}$, by a chain of natural pointed homotopy equivalences. Since
$V\in\cI^{(k)}\subset\cO k$, we have
$$\eta_k\co F(V)@>\simeq>>T_kF(V)$$
so that $L_kF(V)$ is homotopy equivalent to the homotopy
fiber of the map $\eta_j\co F(V)\la T_jF(V)$. Recall that $T_jF(V)$
is defined as 
$$\holim_{U\in\cO j(V)}F(U)\,.$$
Now observe that the inclusion of posets
$$\{\cup_{s\in R}V_s\mid R\subset\pi_0(V), R\ne\pi_0(V)\}
\quad\hookrightarrow\quad\cO j(V)$$
is right cofinal. Complete the proof by applying the
cofinality theorem for homotopy inverse limits.
\qed
\enddemo

In the case of an embedding cofunctor, $F(V)=\emb(V,N)$ as in 1.3,
proposition 9.1 can be made much more explicit. 
We need a base point in $\emb(M,N)$, so we may as well assume that
$M$ is a smooth submanifold of $N$.  
For $S\in\binom Mk$ let $\Psi(S)$ be the total homotopy fiber of the 
$k$--cube of pointed Spaces
$$\bigl\{\emb(R,N)\mid R\subset S\bigr\}\,.\tag{$**$}$$
These Spaces are pointed because $R\subset S\subset M\subset N$. 

\proclaim{9.2 Theorem} {\sl For $k\ge2$, the homogeneous cofunctor 
$L_k\emb(\bla,N)$ is classified by the fibration $p\co Z\to\binom Mk$ with 
fibers 
$p^{-1}(S)=\Psi(S)\,$.}
\endproclaim

\demo{{\bf Proof}} The first and most important observation here is that,
for every $V$ in $\cI^{(k)}$ and every $S\in\binom Mk$ which has exactly
one point in each component of $V$, the obvious restriction map
$$\rho\co \Phi(V)\la\Psi(S)$$
is a homotopy equivalence. This can be seen as follows. For each 
$R\subset S$, there is a homotopy pullback square
$$
\CD
\emb(V_R,N)@>>>\mono(TM|_R,TN) \\
@VVV     @VVV \\
\emb(R,N) @>\subset>> N^R
\endCD$$
where $\mono$ denotes a Space of bundle monomorphisms (vector bundle 
morphisms which are mono in each fiber of the domain). 
Allowing $R$ to be a variable subset of $S$
we may think of it as a square in which each vertex is a $k$--cube 
of Spaces. The total homotopy fibers of these $k$--cubes will then again 
form a homotopy pullback square. But the two $k$--cubes in the right--hand
column are fibrant, so their total homotopy fibers agree with their
total fibers, which reduce to a single point if $k\ge2$ (but not if $k=1$).
Therefore
the total homotopy fibers of the $k$--cubes in the left--hand column are 
homotopy equivalent, which amounts to saying that $\rho\co \Phi(V)\to\Psi(V)$
is a homotopy equivalence.

Now let $E_{\flat}$ be the homotopy colimit of the cofunctor taking 
$V\in\cI^{(k)}$ to the space of sections of $V\to\pi_0(V)$. There are 
obvious forgetful maps
$$|\cI^{(k)}|@<<<E_{\flat}@>>>\tbinom Mk\,.$$
The first of these is a homotopy equivalence by inspection. Comparison 
with the space $E$ in the proof of 3.5 (towards the end) shows that
the second map is also a homotopy equivalence. In more detail,
there is a commutative diagram
$$\CD  
|\cI^{(k)}|@<\simeq<<E_{\flat}@>>>\tbinom Mk \\
@VV=V   @VV\subset V    @VV=V \\
|\cI^{(k)}|@<\simeq<<E@>\simeq>>\tbinom Mk\,.
\endCD$$
Let $p_1$ be the pullback of the quasifibration in 9.1 to $E_{\flat}$,
and let $p_2$ be the pullback of the fibration in 9.2 to $E_{\flat}$.
From the observation made at the beginning of this proof, it is clear
that there is a map over $E_{\flat}$ from $p_1$ to $p_2$ which maps
each fiber of $p_1$ to the corresponding fiber of $p_2$ by a homotopy
equivalence. \qed
\enddemo

In proposition 9.2, the case $k=1$ has been excluded because it is 
different. However, it is also well understood:  We have
$$T_1\emb(\bla,N)=L_1\emb(\bla,N)\,\simeq\,\imm(\bla,N)\,.$$
This follows easily from 5.1 and the observation that all arrows
in the commutative square
$$\CD
\emb(V,N)@>\subset>> \imm(V,N) \\
@VV\eta_1 V          @VV\eta_1 V \\
T_1\emb(V,N)@>\subset>> T_1\imm(V,N)
\endCD$$
are homotopy equivalences if $V\in\cO_1$. 

\section{Boundary Conditions}
So far all manifolds considered were without boundary. When there are 
boundaries, the theory looks slightly different. The following is an 
outline. 

Suppose that $M^m$ is smooth, possibly with boundary. 
Let $\cO$ be the poset of all open subsets of $M$
which contain $\partial M$. 
A cofunctor from $\cO$ to Spaces is {\sl good}
if it satisfies conditions (a) and (b) just before 1.2, literally. 
In (a) we use a definition of {\sl isotopy equivalence} which is 
appropriate for manifolds with boundary:  a smooth codimension 
embedding $(V,\partial V)\to (W,\partial W)$ is an isotopy equivalence 
if, and so on.
 
\definition{10.1 Example} Suppose that $M$ is a neat smooth   
submanifold of another smooth manifold $N$ with boundary. 
That is, $M$ meets $\partial N$ 
transversely, and $\partial M=M\cap\partial N$. For $V$ in $\cO$ let $F(V)$
be the Space of smooth embeddings $V\to N$ which agree with the inclusion
near $\partial M\subset V$. Then $F$ is good. 
\enddefinition

\definition{10.2 Example} Suppose that $M$ is a smooth submanifold with
boundary of another smooth manifold $N$ {\sl without} boundary. For 
$V$ in $\cO$ let $F(V)$
be the Space of smooth embeddings $V\to N$ which agree with the inclusion
near $\partial M\subset V$. Then $F$ is good. 
\enddefinition

In practice example 10.1 is more important because it cannot be reduced 
to simpler cases, whereas 10.2 can often be so
reduced. For example, with $F$ as in 10.2 there is 
a fibration sequence up to homotopy
$$F(M)@>>>\emb(M\minus\partial M,N)@>>>\emb(\partial M,N)$$
provided $\partial M$ is compact. This follows from the 
isotopy extension theorem. It is a mistake to think that a 
similar reduction is possible in the case of 10.1. (Unfortunately 
I made that mistake in \cite{23, \S5}, trying to avoid further
definitions; the calculations done there are nevertheless correct.)

In both examples, 10.1 and 10.2, the values $F(V)$ are contractible for 
collar neighborhoods $V$ of $\partial M$. For general $F$, this may not 
be the case. 

The definition of a {\sl polynomial} cofunctor of some degree $\le k$ is 
again literally the same as before (2.2); we must 
insist that the closed subsets $A_0,\dots A_k$ of $V\in\cO$ have 
empty intersection with $\partial M$, since otherwise 
$F(V\minus\cup_{i\in S}A_i)$ is not defined. 

The definition of the full subcategory $\cO k$ is more complicated.
An element $V\in\cO$ belongs to $\cO k$ if it is a union 
of two disjoint open subsets $V_1$ and $V_2$, where $V_1$ is a collar 
about $\partial M$ (diffeomorphic to $\partial M\times[0,1)$) and $V_2$
is diffeomorphic to a disjoint union of $\le k$ copies of $\RR^m$. 

Later we will need a certain subcategory $\cI^{(k)}$ of $\cO k$.
An object of $\cO k$ belongs to $\cI^{(k)}$ if it has 
exactly $k$ components not meeting $\partial M$; the morphisms in $\cI^{(k)}$ are 
the inclusions which are isotopy equivalences.

As before, $T_kF$ can be defined as the homotopy right Kan extension 
along $\cO k\to\cO$ of $F|\cO k$. It turns out to be polynomial of 
degree $\le k$, and it turns out that $\eta_k\co F\to T_kF$ has the 
properties listed in 6.1. 

If $F(M)$ comes with a selected base point, then we can define $L_kF(V)$
as the homotopy fiber of $T_kF(V)\to T_{k-1}F(V)$. The cofunctor $L_kF$
is homogeneous of degree $k$ (definition like 8.1). 

A general procedure for making homogeneous cofunctors of degree $k$ on
$\cO$ is as follows. {\sl Notation:} $\iota$ is the ``delete boundary'' 
command. Let $p\co Z\to\binom{\iota M}k$ be a fibration. Suppose
that it has a distinguished partial 
section defined near $K$, where $K$ consists of all the points in the 
symmetric product $\spr_kM$ having at least two identical coordinates,
or having at least one coordinate in $\partial M$. For $V$ in $\cO$ let 
$E(V)$ be the Space of (partial) sections of $p$ defined over 
$\binom{\iota V}k$ which agree with the distinguished (zero) section near $K$. 
Then $E$ is homogeneous of degree $k$. 

There is a classification theorem for homogeneous cofunctors 
of degree $k$ on $\cO$, 
to the effect that up to equivalence they can all be obtained in the 
way just described. The classifying fibration $p$ for a homogeneous $E$
of degree $k$ can be found/recovered as follows.
Suppose that $S\subset\iota M$ has $k$ elements.
Choose $V\in \cO k$ so that $V$ contains $S\cup\partial M$ as a 
deformation retract. For $R\subset S$, let $V_R$ be the union of the components
of $V$ which meet $\partial M\cup R$. Let $\Phi(V)$ be the total homotopy 
fiber of the $k$--cube
$$R\mapsto F(V_R)\,.$$
Then $\Phi(V)\simeq p^{-1}(S)$. If more detailed information is needed,
one has to resort to quasifibrations:  the rule $V\mapsto \Phi(V)$ can 
be regarded as a cofunctor on $\cI^{(k)}$
and it gives rise to a quasifibration on $|\cI^{(k)}|\simeq\binom{\iota M}k$. 
The associated fibration is $p$. 

\definition{10.3 Example} In the situation of 10.1, the classifying 
fibration $p_k$ for $L_kF$ has $p_k^{-1}(S)$ equal to the total homotopy 
fiber of the $k$--cube
$$R\mapsto\emb(R,\iota N)$$
for $R\subset S$, provided $k\ge2$. The case $k=0$ is uninteresting
(fiber contractible, base a single point). The case $k=1$ is 
different as usual; for $s\in\iota M$, the fiber 
$p_1^{-1}(\{s\})$ is the space of linear monomorphisms $T_sM\to T(\iota N)$.
All this is exactly as in 9.2. 
For example, suppose that $M$ is compact (with boundary). 
Then $L_kF(M)$ is homotopy 
equivalent to the space of sections of $p_k$ with {\sl compact support}.
In other words, we are dealing with sections defined on all of the configuration space
$\binom{\iota M}k$ and
equal to the zero section outside a compact set. 
\enddefinition

\subheading{Acknowledgments} Tom Goodwillie has exerted a very strong 
influence on this work---not only by kindly communicating his ideas of long 
ago to me, but also by suggesting the right analogies at the right time.
Guowu Meng patiently discussed his calculations of certain embedding
Spaces with me. Finally I am greatly indebted to John Klein for drawing 
my attention to Meng's thesis and Goodwillie's disjunction theory, on which
it is based.  

The author is partially supported by the NSF.


\Refs
\ref
\no1\by A\,K Bousfield\by D\,M Kan
\book Homotopy limits, completions and localizations
\bookinfo Lect. Notes in Math. vol. 304, Springer \yr1973
\endref

\ref 
\no2\by A Dold \book Lectures on Algebraic Topology 
\bookinfo Grundlehren series, \publ Springer \yr1972
\endref

\ref 
\no3\by W Dwyer \paper The centralizer decomposition of BG
\inbook Algebraic topology: new trends in localization and 
periodicity (San Feliu de Guixols, 1994) \bookinfo Progr. in Math. \vol136
\publ Birkh\"auser \yr1996 \pages 167--184
\endref

\ref 
\no4\by W Dwyer \by D\,M Kan \paper A classification theorem for 
diagrams of simplicial sets
\jour Topology, \vol23 \yr1984 \pages139--155
\endref

\ref
\no5\by W Dwyer\by M Weiss \by B Williams \paper A parametrized
index theorem for the algebraic K--theory Euler class, Parts I + II
\paperinfo preprint (1998)
\endref

\ref 
\no6\by T Goodwillie \paper Calculus I: The first derivative of
pseudoisotopy theory \jour K--Theory, \vol4 \yr1990 \pages1--27
\endref

\ref
\no7\by T Goodwillie \paper Calculus II: Analytic Functors
\jour K--Theory, \vol5 \yr1992 \pages295--332
\endref



\ref \no8\by T Goodwillie \paper Calculus III: The Taylor series of a 
homotopy functor \paperinfo in preparation 
\endref

\ref \no9\by T Goodwillie \paper A multiple disjunction lemma for 
smooth concordance embeddings \jour Memoirs Amer. Math. Soc.
\vol 86 no. 431 \yr1990
\endref

\ref \no10\by T Goodwillie \paper Excision estimates for spaces of 
homotopy equivalences \paperinfo pre\-print, Brown University
(1995)
\endref

\ref \no11\by T Goodwillie \paper Excision estimates for spaces of 
diffeomorphisms \paperinfo preprint, Brown University (1998)
\endref

\ref \no12\by T Goodwillie \by J Klein \paper Excision estimates for 
spaces of Poincar\'e embeddings \paperinfo in preparation
\endref

\ref
\no13\by T Goodwillie \by M Weiss \paper Embeddings from the point 
of view of immersion theory: Part II \jour Geometry and Topology,
3
\yr 1999\pages 103--118 
\endref

\ref 
\no14\by M\,L Gromov \paper Stable mappings of foliations into manifolds
\jour Izv. Akad. Nauk. SSSR Mat. \vol33 \yr1969 \pages707--734, \moreref
\jour Trans. Math. USSR (Izvestia) \vol3 \yr1969 \pages671--693
\endref


\ref
\no15\by A Haefliger \by V Poenaru \paper Classification des immersions
combinatoires
\jour Publ.\linebreak Math. I.H.E.S. \vol23  \yr1964 \pages75--91
\endref

\ref
\no16\by R Kirby \by L Siebenmann
\book Foundational Essays on Topological Manifolds, \linebreak Smoothings,
and Triangulations
\bookinfo Ann. of Math. Studies vol. 88, \publ Princeton Univ. Press
\yr1977
\endref

\ref
\no17\by S Mac Lane \book Categories for the Working Mathematician
\publ Springer \yr1971
\endref

\ref
\no18\by S Mac Lane \by I Moerdijk \book Sheaves in Geometry and Logic
\publ Springer--Verlag \yr1992
\endref

\ref
\no19\by V Poenaru \paper Homotopy and differentiable singularities
\inbook Manifolds--Am\-sterdam 1970 
\bookinfo Springer Lect. Notes in
Math. vol. 197, \pages106--132
\endref

\ref \no20\by G Segal \paper Classifying spaces related 
to foliations \jour Topology, \vol13 \yr 1974 \pages293--312
\endref

\ref
\no21\by R Thomason \paper Homotopy colimits in the category of small 
categories \jour Math. Proc. Camb. Phil. Soc. \vol85 \yr1979 \pages91--109
\endref

\ref 
\no22\by F Waldhausen \paper Algebraic K--theory of Spaces
\inbook Proc. of 1983 Rutgers Conf. on Algebraic Topology 
\bookinfo Lect. Notes in Math. \vol1126 \yr1985 \pages318--419
\endref

\ref
\no23\by M Weiss \paper Calculus of Embeddings
\jour Bull. Amer. Math. Soc. \vol33\yr1996 \pages177--187
\endref

\ref
\no24\by M Weiss \paper Curvature and finite domination
\jour Proc. Amer. Math. Soc. \vol124 \yr1996 \pages615--622
\endref

\ref
\no25\by M Weiss \paper Orthogonal Calculus
\jour Trans. Amer. Math. Soc. \vol347\yr1995 \pages3743--3796,
\moreref \paper Erratum \jour Trans. Amer. Math. Soc. \vol350 
\yr1998 \pages851--855 
\endref
\endRefs
\enddocument